\font \sevenrm=cmr7
\font \fiverm=cmr5
\documentclass[11pt, english]{amsart}
\usepackage{amsfonts}
\usepackage{amssymb}
\usepackage{amsmath}
\usepackage[all]{xy}               %xypic macro for latex2.09
\usepackage{color}
\usepackage{graphicx}
\usepackage{xspace}
\usepackage{axodraw}
\def\diagramme #1{\vskip 4mm \centerline {#1} \vskip 4mm}

 \newcommand{\nc}{\newcommand}

 {\everymath{\displaystyle\everymath{}}\array}%
 {\endarray}

\newtheorem{thm}{Theorem}

\nc{\comment}[1]{[[{\tt #1}]] }
\nc{\Cal}[1]{{\mathcal {#1}}}
\nc{\mop}[1]{\mathop{\hbox {\rm #1} }\nolimits}
\nc{\gmop}[1]{\mathop{\hbox {\bf #1} }\nolimits}

\nc{\smop}[1]{\mathop{\hbox {\sevenrm #1} }\nolimits}
\nc{\ssmop}[1]{\mathop{\hbox {\fiverm #1} }\nolimits}
\nc{\mopl}[1]{\mathop{\hbox {\rm #1} }\limits}
\nc{\smopl}[1]{\mathop{\hbox {\sevenrm #1} }\limits}
\nc{\ssmopl}[1]{\mathop{\hbox {\fiverm #1} }\limits}
\nc{\frakg}{{\frak g}}
\nc{\g}[1]{{\frak {#1}}}
\def \restr#1{\mathstrut_{\textstyle |}\raise-6pt\hbox{$\scriptstyle #1$}}
\def \srestr#1{\mathstrut_{\scriptstyle |}\hbox to
  -1.5pt{}\raise-4pt\hbox{$\scriptscriptstyle #1$}}
\nc{\wt}{\widetilde} \nc{\wh}{\widehat}

\nc{\redtext}[1]{\textcolor{red}{\tt [[#1]]}}
\nc{\bluetext}[1]{\textcolor{blue}{#1}}

\nc\fleche[1]{\mathop{\hbox to #1 mm{\rightarrowfill}}\limits}
\nc{\ignore}[1]{}
\def\semi{\mathrel{\times}\kern -.85pt\joinrel\mathrel{\raise
    1.4pt\hbox{${\scriptscriptstyle |}$}}}
\nc\R{{\mathbb R}}
\nc\N{{\mathbb N}}
\nc\inver{^{-1}}
\nc\point{\hbox{\bf .}}
\nc\un{\hbox{\bf 1}}
\def\link#1#2{\raise -2pt\hbox{$\scriptstyle #1-\!\!-\!\!- #2$}}
\def\slink#1#2{\raise -1.5pt\hbox{$\scriptscriptstyle #1-\!\!\!-\!\!\!- #2$}}
\def\HRT{\Cal H^{\Cal D}_{\smop{BCK}}}
\def\HRTX{\Cal H^{X}_{\smop{BCK}}}
\def\HRTY{\Cal H^{Y}_{\smop{BCK}}}
\def\GRT{G^{\Cal D}_{\smop{BCK}}}

\def\racine{{\scalebox{0.3}{ 
\begin{picture}(12,12)(38,-38)
\SetWidth{0.5} \SetColor{Black} \Vertex(45,-33){5.66}
\end{picture}}}}

 \def\arbrea{\,{\scalebox{0.15}{ 
  \begin{picture}(8,55) (370,-248)
    \SetWidth{2}
    \SetColor{Black}
    \Line(374,-244)(374,-200)
    \Vertex(374,-197){9}
    \Vertex(375,-245){12}
  \end{picture}
}}\,}

 \def\arbreba{\,{\scalebox{0.15}{ 
\begin{picture}(8,106) (370,-197)
    \SetWidth{2}
    \SetColor{Black}
    \Line(374,-193)(374,-149)
    \Vertex(374,-146){9}
    \Vertex(375,-194){12}
    \Line(374,-142)(374,-98)
    \Vertex(374,-95){9}
  \end{picture}
}}\,}

 \def\arbrebb{\,{\scalebox{0.15}{ 
  \begin{picture}(48,48) (349,-255)
    \SetWidth{2}
    \SetColor{Black}
    \Vertex(375,-252){12}
    \Line(376,-250)(395,-215)
    \Line(373,-251)(354,-214)
    \Vertex(353,-211){9}
    \Vertex(395,-213){9}
  \end{picture}
}}}

\def\arbrebbdec#1#2#3{\,{\scalebox{0.60}{ 
  \begin{picture}(48,48) (349,-255)
    \SetWidth{1}
    \SetColor{Black}
    \Vertex(375,-252){12}
    \Line(376,-250)(395,-215)
    \Line(373,-251)(354,-214)
    \Vertex(353,-213){9}
    \Vertex(395,-213){9}
    \SetColor{White}
    \Vertex(375,-252){11}
    \Vertex(353,-213){8}
    \Vertex(395,-213){8}
    \Text(371,-257)[lb]{\Large{\Black{$#1$}}}
    \Text(390,-218)[lb]{\Large{\Black{$#3$}}}
    \Text(348,-218)[lb]{\Large{\Black{$#2$}}}
  \end{picture}
}}\ }

\def\arbreca{\,{\scalebox{0.15}{
\begin{picture}(8,156) (370,-147)
    \SetWidth{2}
    \SetColor{Black}
    \Line(374,-143)(374,-99)
    \Vertex(374,-96){9}
    \Vertex(375,-144){12}
    \Line(374,-92)(374,-48)
    \Vertex(374,-45){9}
    \Line(374,-42)(374,2)
    \Vertex(374,5){9}
  \end{picture}
}}\,}

\def\arbrecb{\,{\scalebox{0.15}{
\begin{picture}(48,94) (349,-255)
\SetWidth{2}
\SetColor{Black}
\Line(376,-204)(395,-169)
\Line(373,-205)(354,-168)
\Vertex(353,-165){9}
\Vertex(395,-167){9}
\Vertex(374,-205){9}
\Line(374,-246)(374,-209)
\Vertex(374,-252){12}
\end{picture}}}\,}

\def\arbrecc{\,{\scalebox{0.15}{
 \begin{picture}(48,98) (349,-205)
    \SetWidth{2}
    \SetColor{Black}
    \Vertex(375,-202){12}
    \Line(376,-200)(395,-165)
    \Line(373,-201)(354,-164)
    \Vertex(353,-161){9}
    \Vertex(395,-163){9}
    \Line(353,-160)(353,-113)
    \Vertex(353,-111){9}
  \end{picture}
}}\,}

\def\arbreccbis{\,{\scalebox{0.15}{
 \begin{picture}(48,98) (349,-205)
    \SetWidth{2}
    \SetColor{Black}
    \Vertex(375,-202){12}
    \Line(376,-200)(395,-165)
    \Line(373,-201)(354,-164)
    \Vertex(353,-161){9}
    \Vertex(395,-163){9}
    \Line(353,-160)(353,-113)
    \Vertex(353,-111){9}
    \SetColor{White}
    \Vertex(375,-202){10}
    \Vertex(353,-161){7}
  \end{picture}
}}\,}

\def\arbreccter{\,{\scalebox{0.15}{
 \begin{picture}(48,98) (349,-205)
    \SetWidth{2}
    \SetColor{Black}
    \Vertex(375,-202){12}
    \Line(376,-200)(395,-165)
    \Line(373,-201)(354,-164)
    \Vertex(353,-161){9}
    \Vertex(395,-163){9}
    \Line(353,-160)(353,-113)
    \Vertex(353,-111){9}
    \SetColor{White}
     \Vertex(375,-202){10}
  \end{picture}
}}\,}

\def\arbrecd{\,{\scalebox{0.15}{
\begin{picture}(48,52) (349,-251)
    \SetWidth{2}
    \SetColor{Black}
    \Vertex(375,-248){12}
    \Line(376,-246)(395,-211)
    \Line(373,-247)(354,-210)
    \Vertex(353,-207){9}
    \Vertex(395,-209){9}
    \Line(375,-247)(375,-206)
    \Vertex(376,-203){9}
  \end{picture}
 }}\,}

\def\arbreda{\,{\scalebox{0.15}{
\begin{picture}(8,204) (370,-99)
    \SetWidth{2}
    \SetColor{Black}
    \Line(374,-95)(374,-51)
    \Vertex(374,-48){9}
    \Vertex(375,-96){12}
    \Line(374,-44)(374,0)
    \Vertex(374,3){9}
    \Line(374,6)(374,50)
    \Vertex(374,53){9}
    \Line(374,53)(374,98)
    \Vertex(374,101){9}
  \end{picture}
}}\,}

\def\arbredb{\,{\scalebox{0.15}{
\begin{picture}(48,135) (349,-255)
    \SetWidth{2}
    \SetColor{Black}
    \Line(376,-163)(395,-128)
    \Line(373,-164)(354,-127)
    \Vertex(353,-124){9}
    \Vertex(395,-126){9}
    \Vertex(374,-164){9}
    \Line(374,-205)(374,-168)
    \Vertex(374,-207){9}
    \Line(374,-248)(374,-211)
    \Vertex(374,-252){12}
  \end{picture}
}}\,}

\def\arbredc{\,{\scalebox{0.15}{
 \begin{picture}(48,150) (349,-205)
    \SetWidth{2}
    \SetColor{Black}
    \Line(376,-148)(395,-113)
    \Line(373,-149)(354,-112)
    \Vertex(353,-109){9}
    \Vertex(395,-111){9}
    \Line(353,-108)(353,-61)
    \Vertex(353,-59){9}
    \Line(374,-200)(374,-153)
    \Vertex(374,-149){9}
    \Vertex(374,-202){12}
  \end{picture}
}}\,}

\def\arbredd{\,{\scalebox{0.15}{
 \begin{picture}(48,99) (349,-251)
    \SetWidth{2}
    \SetColor{Black}
    \Line(376,-199)(395,-164)
    \Line(373,-200)(354,-163)
    \Vertex(353,-160){9}
    \Vertex(395,-162){9}
    \Vertex(376,-156){9}
    \Vertex(376,-248){12}
    \Line(375,-245)(375,-204)
    \Line(375,-200)(375,-159)
    \Vertex(375,-201){9}
  \end{picture}
}}\,}

\def\arbrede{\,{\scalebox{0.15}{
 \begin{picture}(48,153) (349,-150)
    \SetWidth{2}
    \SetColor{Black}
    \Vertex(375,-147){12}
    \Line(376,-145)(395,-110)
    \Line(373,-146)(354,-109)
    \Vertex(353,-106){9}
    \Vertex(395,-108){9}
    \Line(353,-105)(353,-58)
    \Vertex(353,-56){9}
    \Line(353,-52)(353,-5)
    \Vertex(353,-1){9}
  \end{picture}
}}\,}

\def\arbredf{\,{\scalebox{0.15}{
\begin{picture}(48,98) (349,-205)
    \SetWidth{2}
    \SetColor{Black}
    \Vertex(375,-202){12}
    \Line(376,-200)(395,-165)
    \Line(373,-201)(354,-164)
    \Vertex(353,-161){9}
    \Vertex(395,-163){9}
    \Line(353,-160)(353,-113)
    \Vertex(353,-111){9}
    \Line(395,-159)(395,-112)
    \Vertex(395,-111){9}
  \end{picture}
}}\,}

\def\arbredz{\,{\scalebox{0.15}{
  \begin{picture}(68,88) (329,-215)
    \SetWidth{2}
    \SetColor{Black}
    \Vertex(375,-212){12}
    \Line(376,-210)(395,-175)
    \Line(373,-211)(354,-174)
    \Vertex(353,-171){9}
    \Vertex(395,-173){9}
    \Line(351,-168)(332,-131)
    \Line(355,-168)(374,-133)
    \Vertex(333,-131){9}
    \Vertex(374,-131){9}
  \end{picture}
}}\,}

\def\arbredg{\,{\scalebox{0.15}{
\begin{picture}(48,98) (349,-205)
    \SetWidth{2}
    \SetColor{Black}
    \Vertex(375,-202){12}
    \Line(376,-200)(395,-165)
    \Line(373,-201)(354,-164)
    \Vertex(353,-161){9}
    \Vertex(395,-163){9}
    \Line(375,-201)(375,-160)
    \Vertex(376,-157){9}
    \Vertex(376,-111){9}
    \Line(375,-155)(375,-114)
  \end{picture}
}}\,}

\def\arbredh{\,{\scalebox{0.15}{
 \begin{picture}(90,46) (330,-257)
    \SetWidth{2}
    \SetColor{Black}
    \Vertex(375,-254){12}
    \Line(376,-252)(395,-217)
    \Vertex(395,-215){9}
    \Line(374,-254)(335,-226)
    \Vertex(334,-224){9}
    \Line(375,-252)(356,-215)
    \Vertex(355,-215){9}
    \Line(374,-255)(417,-227)
    \Vertex(418,-225){9}
  \end{picture}
}}\,}

\def\shu{\joinrel{\!\scriptstyle\amalg\hskip -3.1pt\amalg}\,}
\def\qshu{\joinrel{\!\scriptstyle\amalg\hskip -3.1pt\amalg}\,\hskip -8pt\hbox{-}\hskip 5pt}
\def\sshu{\joinrel{\,\scriptscriptstyle\amalg\hskip -2.5pt\amalg}\,}
\def\sqshu{\joinrel{\,\scriptscriptstyle\amalg\hskip -2.5pt\amalg}\,\hskip -7pt\hbox{-}\hskip 4pt}
\def\sqshubis{\joinrel{\,\scriptscriptstyle\amalg\hskip -2.5pt\amalg}\,\hskip -8pt\hbox{-}\hskip 4pt}
\def\surj{\to\hskip -9pt\to}
\def\ct{T\hskip -5pt\hbox{-}}
  
\begin{document}
%%%%%%%%%%%%%%%%%%%%%%%%%%%%%%%%%%%%%%%%%%%%%
%%%%%%%%%%%%%%%%%%%%%%%%%%%%%%%%%%%%%%%%%%%%%

\title[Arborified multiple zeta values]{Arborified multiple zeta values}

\author{Dominique Manchon}
\address{Universit\'e Blaise Pascal,
         C.N.R.S.-UMR 6620, 
         3 place Vasar\'ely, BP 80026,
         63178 Aubi\`ere, France}       
         \email{manchon@math.univ-bpclermont.fr}
         \urladdr{http://math.univ-bpclermont.fr/~manchon/}

%%%%%%%%%%%%%%%%%%%%%%%%%%%%%%%%%%%%%%%%%%%%%%%%%%%%%%%%%%%%%%%%%%%
\date{January 25th 2013}
%%%%%%%%%%%%%%%%%%%%%%%%%%%%%%%%%%%%%%%%%%%%%%%%%%%%%%%%%%%%%%%%%%%

\begin{abstract}

We describe some particular finite sums of multiple zeta values which 
arise from J. Ecalle's "arborification", a process which can be
described as a surjective Hopf algebra 
morphism from the Hopf algebra of decorated rooted forests onto a Hopf 
algebra of shuffles or quasi-shuffles. This formalism holds for both the 
iterated sum picture and the iterated integral picture. It involves a 
decoration of the forests by the positive integers in the first case, by 
only two colours in the second case.

\bigskip


\smallskip

\end{abstract}

%
%\begin{altabstract}
%\end{altabstract}
%

\maketitle

%%%%%%%%%%%%%%%%%%%%%%%%%%%%%%%%%%%%%%%%%%%%%%%%%%%%%%%%%%%%%%%%%%%

% \tableofcontents

%\newpage

%%%%%%%%%%%%%%%%%%%%%%%%%%%%%%%%%%%%%%%%%%%%%%%%%%%%%%%%%%%%%%%%%%%

\section{Introduction}
\label{sect:intro}
%%%%%%%%%
\noindent
Multiple zeta values are defined by the following nested sums:
\begin{equation}\label{mz}
\zeta(n_1,\ldots ,n_r):=\sum_{k_1>k_2>\cdots >k_r\ge 1}\frac{1}{k_1^{n_1}\cdots k_r^{n_r}},
\end{equation}
where the $n_j$'s are positive integers. The nested sum \eqref{mz} converges as long as $n_1\ge 2$. The integer $r$ is the \textsl{depth}, whereas the sum $p:=n_1+\cdots n_r$ is the \textsl{weight}. Although the multiple zeta values of depth one and two were already known by L. Euler, the  full set of multiple zeta values first appears in 1981 in a preprint of Jean Ecalle under the name "moule $\zeta_<^\bullet$", in the context of resurgence theory in complex analysis \cite[Page 429]{E2}, together with its companion $\zeta_{\le}^\bullet$ now known as the set of multiple star zeta values. The systematic study begins a decade later with the works of M.~E.~Hoffman \cite{H} and D.~Zagier \cite{Z}. It has been remarked by M. Kontsevich (\cite{Z}, see also the intriguing precursory Remark 4 on Page 431 in \cite{E2}) that multiple zeta values admit another representation by iterated integrals, namely:
\begin{equation}\label{mz-integral}
\zeta(n_1,\ldots ,n_r)=\int\cdots\int_{0\le u_p\le\cdots u_1\le 1}\,\frac{du_1}{\varphi_1(u_1)}\cdots\frac{du_p}{\varphi_p(u_p)},
\end{equation}
with $\varphi_j(u)=1-u$ if $j\in\{n_1,n_1+n_2,n_1+n_2+n_3,\ldots ,p\}$ and $\varphi_j(u)=u$ otherwise. For later use we set:
$$f_0(u):=u, \hskip12mm f_1(u):=1-u.$$
Iterated integral representation \eqref{mz-integral} is the starting point to the modern approach in terms of mixed Tate motives over ${\mathbb Z}$, already outlined in \cite{Z} and widely developed in the literature since then \cite{T, DG, B1, B2, B3}. Multiple zeta values verify a lot of polynomial relations with integer coefficients: the representation \eqref{mz} by nested sums leads to \textsl{quasi-shuffle relations}, whereas representation \eqref{mz-integral} by iterated integrals leads to \textsl{shuffle relations}. A third family of relations, the \textsl{regularization relations}, comes from a subtle interplay between the two first groups of relations, involving divergent multiple zeta sums $\zeta(1,n_2\ldots n_r)$. A representative example of each family (in the order above) is given by:
\begin{eqnarray}
\zeta(2,3)+\zeta(3,2)+\zeta(5)&=&\zeta(2)\zeta(3),\label{ex-qshu}\\
\zeta(2,3)+3\zeta(3,2)+6\zeta(4,1)&=&\zeta(2)\zeta(3),\label{ex-shu}\\
\zeta(2,1)&=&\zeta(3),\label{euler}
\end{eqnarray}
It is conjectured that these three families include all possible polynomial relations between multiple zeta values. Note that the rationality of the quotient $\displaystyle\frac{\zeta(2k)}{\pi^{2k}}$, proved by L. Euler, does not yield supplementary polynomial identities. As an example, $\zeta(2)=\displaystyle\frac{\pi^2}{6}$ and $\zeta(4)=\displaystyle\frac{\pi^4}{90}$ yield $2\zeta(2)^2=5\zeta(4)$, a relation which can also be deduced from quasi-shuffle, shuffle and regularization relations.\\

It is convenient to write multiple zeta values in terms of words. In view of representations \eqref{mz} and \eqref{mz-integral}, this can be done in two different ways. We consider the two alphabets:
\begin{equation*}
X:=\{x_0,x_1\},\hskip 12mm Y:=\{y_1,y_2,y_3,\ldots\},
\end{equation*}
and we denote by $X^*$ (resp. $Y^*$) the set of words with letters in $X$ (resp. $Y$). The vector space ${\mathbb Q}\langle X\rangle$ freely generated by $X^*$ is a commutative algebra for the \textsl{shuffle product}, which is defined by:
\begin{equation}\label{shu-prod}
(v_1\cdots v_p)\shu(v_{p+1}\cdots v_{p+q}):=\sum_{\sigma\in\smop{Sh}(p,q)}v_{\sigma_1^{-1}}\cdots v_{\sigma_{p+q}^{-1}}
\end{equation}
with $v_j\in X$, $j\in\{1,\ldots ,p+q\}$. Here, $\mop{Sh}(p,q)$ is the set of \textsl{$(p,q$)-shuffles}, i.e. permutations $\sigma$ of $\{1,\ldots p+q\}$ such that $\sigma_1<\cdots \sigma_p$ and $\sigma_{p+1}<\cdots <\sigma_{p+q}$. The vector space ${\mathbb Q}\langle Y\rangle$ freely generated by $Y^*$ is a commutative algebra for the \textsl{quasi-shuffle product}, which is defined as follows: a \textsl{$(p,q)$-quasi-shuffle of type $r$} is a surjection $\sigma:\{1,\ldots p+q\}\surj\{1,\ldots p+q-r\}$ such that  $\sigma_1<\cdots \sigma_p$ and $\sigma_{p+1}<\cdots <\sigma_{p+q}$. Denoting by $\mop{Qsh}(p,q;r)$ the set of $(p,q)$-quasi-shuffles of type $r$, the formula for the quasi-shuffle product $\qshu$ is:
\begin{equation}\label{qshu-prod}
(w_1\cdots w_p)\qshu(w_{p+1}\cdots w_{p+q}):=\sum_{r\ge 0}\,\sum_{\sigma\in\smop{Qsh}(p,q;r)}
w_1^\sigma\cdots w_{p+q-r}^\sigma
\end{equation}
with $w_j\in Y$, $j\in\{1,\ldots ,p+q\}$, and where $w_j^{\sigma}$ is the \textsl{internal product} of the letters in the set $\sigma^{-1}(\{j\})$, which contains one or two elements. The internal product is defined by $[y_ky_l]:=y_{k+l}$.\\

We denote by $Y^*_{\smop{conv}}$ the submonoid of words $w=w_1\cdots w_r$ with $w_1\not =y_1$, and we set $X^*_{\smop{conv}}=x_0X^*x_1$. An injective monoid morphism is given by changing letter $y_n$ into the word $x_0^{n-1}x_1$, namely:
\begin{eqnarray*}
{\frak s}:Y^*&\longrightarrow & X^*\\
y_{n_1}\cdots y_{n_r}&\longmapsto & x_0^{n_1-1}x_1\cdots  x_0^{n_r-1}x_1,
\end{eqnarray*}
and restricts to a monoid isomorphism from $Y^*_{\smop{conv}}$ onto $X^*_{\smop{conv}}$. As notation suggests,  $Y^*_{\smop{conv}}$ and $X^*_{\smop{conv}}$ are two convenient ways to symbolize convergent multiple zeta values through representations \eqref{mz} and \eqref{mz-integral} respectively. The following notation is commonly adopted:
\begin{equation}
\zeta_{\sqshu}(y_{n_1}\cdots y_{n_r}):=\zeta(n_1,\ldots n_r)=:\zeta_{\sshu}\big({\frak s}(y_{n_1}\cdots y_{n_r})\big),
\end{equation}
and extended to finite linear combinations of convergent words by linearity. The relation:
 $$\zeta_{\sqshu}=\zeta_{\sshu}\circ{\frak s}$$
is obviously verified. The quasi-shuffle relations then write:
\begin{equation}\label{qsh-rel}
\zeta_{\sqshu}(w\qshu w')=\zeta_{\sqshu}(w)\zeta_{\sqshu}(w')
\end{equation}
for any $w,w'\in Y^*_{\smop{conv}}$, whereas the shuffle relations write:
\begin{equation}\label{sh-rel}
\zeta_{\sshu}(v\shu v')=\zeta_{\sshu}(v)\zeta_{\sshu}(v')
\end{equation}
for any $v,v'\in X^*_{\smop{conv}}$. By fixing an arbitrary value $\theta$ to $\zeta(1)$ and setting $\zeta_{\sqshu}(y_1)=\zeta_{\sshu}(x_1)=\theta$, it is possible to extend $\zeta_{\sqshu}$, resp. $\zeta_{\sshu}$, to all words in $Y^*$, resp. to $X^*x_1$, such that \eqref{qsh-rel}, resp. \eqref{sh-rel}, still holds. It is also possible to extend $\zeta_{\sshu}$ to a map defined on $X^*$ by fixing an arbitrary value $\theta'$ to $\zeta_{\sshu}(x_0)$, such that \eqref{sh-rel} is still valid. We will stick to $\theta'=\theta$ for symmetry reasons, reflecting the following formal equality between two infinite quantities:
$$
\int_0^1\frac{dt}{t}=\int_0^1\frac{dt}{1-t}.
$$
It is easy to show that for any word $v\in X^*$ or $w\in Y^*$, the expressions $\zeta_{\sshu}(v)$ and $\zeta_{\sqshu}(w)$ are polynomial with respect to $\theta$. It is no longer true that extended $\zeta_{\sqshu}$ coincides with extended $\zeta_{\sshu}\circ{\frak s}$, but the defect can be explicitly written:
\begin{thm}[L. Boutet de Monvel, D. Zagier \cite{Z}]\label{bmz}
There exists an infinite-order inversible differential operator $\rho:\R[\theta]\to\R[\theta]$ such that
\begin{equation}
\zeta_{\sshu}\circ{\frak s}=\rho\circ\zeta_{\sqshubis}.
\end{equation}
The operator $\rho$ is explicitly given by the series:
\begin{equation}
\rho=\exp\left(\sum_{n\ge 2}\frac{(-1)^n\zeta(n)}{n}\left(\frac{d}{d\theta}\right)^n\right).
\end{equation}
\end{thm}
\noindent
In particular, $\rho(1)=1$, $\rho(\theta)=\theta$, and more generally $\rho(P)-P$ is a polynomial of degree $\le d-2$ if $P$ is of degree $d$, hence $\rho$ is inversible. A proof of Theorem \ref{bmz} can be read in numerous references, e.g. \cite{C,IKZ,R}. Any word $w\in Y^*_{\smop{conv}}$ gives rise to Hoffman's regularization relation:
\begin{equation}\label{reg}
\zeta_{\sshu}\big(x_1\shu {\frak s}(w)-{\frak s}(y_1\qshu w)\big)=0,
\end{equation}
which is a direct consequence of Theorem \ref{bmz}. The linear combination of words involved above is convergent, hence \eqref{reg} is a relation between convergent multiple zeta values, although divergent ones have been used to establish it. The simplest regularization relation \eqref{euler} is nothing but \eqref{reg} applied to the word $w=y_2$.\\

Rooted trees can enrich the picture in two ways: first of all, considering a rooted tree $t$ with set of vertices $\Cal V(t)$ and decoration $n_v\in{\mathbb Z}_{>0},\,v\in \Cal V(t)$, we define the associated \textsl{contracted arborified multiple zeta value} by:
\begin{equation}\label{treemz}
\zeta^{\ct}(t):=\sum_{k\in D_t}\ \prod_{v\in\Cal V(t)}\frac{1}{k_v^{n_v}},
\end{equation}
where $D_t$ is made of those maps $v\mapsto k_v\in{\mathbb Z}_{>0}$ such that $k_v<k_w$ if and only if there is a path from the root to $w$ through $v$. The sum \eqref{treemz} is convergent as long as $n_v\ge 2$ if $v$ is a leaf of $t$. The definition is multiplicatively extended to rooted forests. A similar definition can be introduced starting from the integral representation \eqref{mz-integral}: considering a rooted tree $\tau$ with set of vertices $\Cal V(\tau)$ and decoration $e_v\in\{0,1\},\,v\in \Cal V(\tau)$, we define the associated \textsl{arborified multiple zeta value} by:
\begin{equation}\label{treemz-integral}
\zeta^T(\tau):=\int_{u\in \Delta_\tau}\ \prod_{v\in\Cal V(\tau)}\frac{du_v}{f_{e_v}},
\end{equation}
where $\Delta_\tau\subset [0,1]^{|\Cal V(\tau)|}$ is made of those maps $v\mapsto u_v\in [0,1]$ such that $u_v\le u_w$ if and only if there is a path from the root to $w$ through $v$. The integral \eqref{treemz-integral} is convergent as long as $e_v=1$ if $v$ is the root of $\tau$ and $e_v=0$ if $v$ is a leaf of $\tau$. A multiplicative extension to two-coloured rooted forests will also be considered. Arborified multiple zeta values, in this non-contracted form, appear in a recent paper by S. Yamamoto \cite{Y14}.\\

Arborified and contracted arborified multiple zeta values are finite linear combinations of ordinary ones. For example we have~:
$$\zeta^{\ct}(\arbrebbdec{n_3}{n_1}{n_2})=\zeta(n_1,n_2,n_3)+\zeta(n_2,n_1,n_3)+\zeta(n_1+n_2,n_3)$$
and, choosing black for colour $0$ and white for colour $1$:
\begin{eqnarray*}
\zeta^T\left(\arbreccbis\right)&=&2\zeta(3,1)+\zeta(2,2),\\
\zeta^T(\arbreccter)&=&3\zeta(4).
\end{eqnarray*}
The terminology comes from J. Ecalle's \textsl{arborification}, a transformation which admits a "simple" and a "contracting" version \cite{E92, EV}. This transformation is best understood in terms of a canonical surjective morphism from Butcher-Connes-Kreimer Hopf algebra of rooted forests onto a corresponding shuffle Hopf algebra (quasi-shuffle Hopf algebra for the contracting arborification) \cite{FM}.\\

The paper is organized as follows: after a reminder on shuffle and quasi-shuffle Hopf algebras, we describe the two versions of arborification in some detail, and we describe a possible transformation from contracted arborified to arborified multiple zeta values, which can be seen as an arborified version of the map $\frak s$ from words in $Y^*$ into words in $X^*$. A more natural version of this arborified $\frak s$ with respect to the tree structures is still to be found.\\

\noindent
\textbf{Acknowledgements:} Research partly supported by Agence Nationale de la Recherche, project "Carma" ANR-12-BS01-0017. I thank Fr\'ed\'eric Fauvet, Fr\'ed\'eric Menous and Emmanuel Vieillard-Baron for important discussions on arborification.
%%%%%
\section{Shuffle and quasi-shuffle Hopf algebras}
%%%%%
Let $V$ be any commutative algebra on a base field $k$ of characteristic zero. The product on $V$ will be denoted by $(a,b)\mapsto [ab]$. This algebra is not supposed to be unital: in particular any vector space can be considered as a commutative algebra with trivial product $(a,b)\mapsto [ab]=0$. The associated \textsl{quasi-shuffle Hopf algebra} is $\big(T(V),\qshu,\Delta\big)$, where $\big(T(V),\Delta\big)$ is the tensor coalgebra:
$$T(V)=\bigoplus_{k\ge 0}V^{\otimes k}.$$
The indecomposable elements of $V^{\otimes k}$ will be denoted by $v_1\cdots v_k$ with $v_j\in V$. The coproduct $\Delta$ is the deconcatenation coproduct:
\begin{equation}
\Delta(v_1\cdots v_k):=\sum_{r=0}^k v_1\cdots v_r\otimes v_{r+1}\cdots v_k.
\end{equation}
The quasi-shuffle product $\qshu$ is given for any $v_1,\ldots v_{p+q}$ by:
\begin{equation}\label{qshu-prod2}
(v_1\cdots v_p)\qshu(v_{p+1}\cdots v_{p+q}):=\sum_{r\ge 0}\,\sum_{\sigma\in\smop{Qsh}(p,q;r)}
v_1^\sigma\cdots v_{p+q-r}^\sigma
\end{equation}
with $v_j\in Y$, $j\in\{1,\ldots ,p+q\}$, and where $v_j^{\sigma}$ is the internal product of the letters in the set $\sigma^{-1}(\{j\})$, which contains one or two elements. Note that if the internal product vanishes, only ordinary shuffles (i.e. quasi-shuffles of type $r=0$) do contribute to the quasi-shuffle product, which specializes to the shuffle product $\shu$ in this case. The tensor coalgebra endowed with the quasi-shuffle product $\qshu$ is a Hopf algebra which, remarkably enough, does not depend on the particular choice of the internal product \cite{H2}. An explicit Hopf algebra isomorphism $\exp$ from $\big(T(V),\qshu,\Delta\big)$ onto $\big(T(V),\shu,\Delta\big)$ is given in \cite{H2}. Although we won't use it, let us recall its expression: let ${\Cal P}(k)$ be the set of compositions of the integer $k$, i.e. the set of
sequences  $I=(i_1,\ldots,i_r)$ of positive integers such that
$i_1+\cdots +i_r=k$. For any  $u=v_1\ldots v_k\in T(V)$
and any composition  $I=(i_1,\ldots,i_r)$ of $r$ we set:
$$I[u]:=[v_1\ldots v_{i_1}].
[v_{i_1+1}\cdots v_{i_1+i_2}]\ldots
[v_{i_1+\cdots+i_{r-1}+1}\ldots v_k].$$
Then:
$$\exp u=\sum_{I=(i_1,\ldots ,i_r)\in{\Cal P}(k)}\frac{1}{i_1!\ldots
  i_r!}I[u].$$
Moreover (\cite{H2}, lemma 2.4), the inverse $\log$ of $\exp$ is given by~:
$$\log u=\sum_{I=(i_1,\ldots ,i_r)\in{\Cal P}(k)}\frac{(-1)^{k-r}}{i_1\ldots
  i_r}I[u].$$
For example for  $v_1,v_2,v_3\in V$ we have:
\begin{eqnarray*}
\exp v_1 =v_1 &,& \log v_1 =v_1,\\
\exp (v_1v_2) =v_1v_2+\frac{1}{2}[v_1 v_2]&,&
\log (v_1 v_2) =v_1 v_2-\frac{1}{2}[v_1v_2],\\
\exp (v_1v_2v_3) &=&v_1v_2v_3+\frac{1}{2}([v_1v_2]
v_3+v_1[v_2v_3])+\frac{1}{6}[v_1 v_2 v_3],\\
\log (v_1 v_2 v_3) &=&v_1 v_2 v_3-\frac{1}{2}([v_1v_2]
v_3+v_1[ v_2v_3])+\frac{1}{3}[v_1 v_2v_3].
\end{eqnarray*}
Going back to the notations of the introduction, ${\mathbb Q}\langle Y \rangle$ is the quasi-shuffle Hopf algebra associated to the algebra $tk[t]$ of polynomials without constant terms, whereas  ${\mathbb Q}\langle X \rangle$ is the shuffle Hopf algebra associated with the two-dimensional vector space spanned by $X$.
%%%%%
\section{The Butcher-Connes-Kreimer Hopf algebra of decorated rooted trees}
%%%%%
Let $\Cal D$ be a set. A \textsl{rooted tree} is an oriented (non planar) graph with a finite number 
of vertices, among which one is distinguished and called the \textsl{root}, 
such that any vertex admits exactly one incoming edge, except the root which 
has no incoming edges. A \textsl{$\Cal D$-decorated rooted tree} is a rooted tree $t$ together with a map from its set of vertices $\Cal V(t)$ into $\Cal D$. Here is the list of (non-decorated) rooted trees up to five vertices:
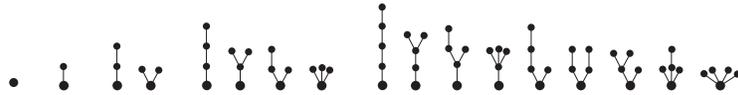
\begin{equation*}
\racine \hskip 5mm \arbrea \hskip 5mm  \arbreba\  \arbrebb \hskip 5mm  
\arbreca\  \arbrecb\  \arbrecc\ \arbrecd \hskip 5mm  \arbreda\  \arbredb\  
\arbredc\  \arbredd\  \arbrede\  \arbredf\  \arbredz\  \arbredg\ \arbredh
\end{equation*}
A \textsl{$\Cal D$-decorated rooted forest} is a finite collection of $\Cal D$-decorated rooted trees, with possible repetitions. 
The empty set is the forest containing no trees, and is denoted by $\un$. 
For any $d\in\Cal D$, the \textsl{grafting operator} $B_+^d$ takes any forest and changes it into a 
tree by grafting all components onto a common root decorated by $d$, with the convention 
$B_+^d(\un)=\racine_d$.\\

Let $\Cal T^{\Cal D}$ denote the set of nonempty rooted trees and let $\HRT=k[\Cal T^{\Cal D}]$ 
be the free commutative unital algebra generated by elements of $\Cal T^{\Cal D}$. 
We identify a product of trees with the forest containing these trees. 
Therefore the vector space underlying $\HRT$ is the linear span of rooted 
forests. This algebra is a graded and connected Hopf algebra, called the 
\textsl{Hopf algebra of $\Cal D$-decorated rooted trees}, with the following structure:
the grading is given by the number of vertices, and
the coproduct on a rooted forest $u$ is 
described as follows \cite{F, Mu}: the set $\Cal V(u)$ of vertices of a forest $u$ is 
endowed with a partial order defined by $x \le y$ if and only if there is a 
path from a root to $y$ passing through $x$. Any subset $W$ of  $\Cal V(u)$ defines a {\sl subforest\/} $u\restr{W}$ of $u$ 
in an obvious manner, i.e. by keeping the edges of $u$ which link two elements 
of $W$. The coproduct is then defined by:
\begin{equation}
\label{coprod}
\Delta(u)= \sum_{V \amalg W=\Cal V(u) \atop W<V} u\restr{V}\otimes u\restr{W}.
\end{equation}
Here the notation $W<V$ means that $y<x$ for any vertex $x$ in $V$ and any
vertex $y$ in $W$ such that $x$ and $y$ are comparable. Such a couple $(V,W)$
is also called an \textsl{admissible cut}, with \textsl{crown} (or pruning) 
$u\restr{V}$ and \textsl{trunk} $u\restr{W}$. We have for example:
 \allowdisplaybreaks{
\begin{eqnarray*}
\Delta\big(\arbrea\big) &=& 
   \arbrea \otimes \un + \un \otimes \arbrea + \racine \otimes \racine \\
\Delta\big(\! \arbrebb \big) &=& \arbrebb \otimes \un + \un \otimes \arbrebb + 
   2\racine \otimes\arbrea + \racine\racine\otimes \racine . 
\end{eqnarray*}}
The counit is $\varepsilon(\un)=1$ and $\varepsilon(u)=0$ for any non-empty 
forest $u$. The coassociativity of the coproduct is easily checked using the following
formula for the iterated coproduct :
$$
\wt\Delta^{n-1}(u)=
\sum_{V_1\amalg\cdots\amalg V_n=\Cal V(u) \atop 
V_n<\cdots <V_1} 
u\restr{V_1}\otimes\cdots\otimes u\restr{V_n}.
$$ 
The notation $V_n<\cdots <V_1$ is to be understood as $V_i<V_j$ for 
any $i>j$, with $i,j\in\{1,\ldots ,n\}$.\\

This Hopf algebra first appeared in the work of A. D\" ur in 1986 \cite{Dur}. It has been rediscovered and intensively studied by D. Kreimer in 1998 
\cite{Kreimer}, as the Hopf algebra describing the combinatorial part 
of the BPHZ renormalization procedure of Feynman graphs in a scalar 
$\varphi^3$ quantum field theory.  Its group of characters:
\begin{equation}
\GRT=\mop{Hom}_{\smop{alg}}(\HRT, k)
\end{equation}
is known as the \textsl{Butcher group} and plays a key role in approximation methods in numerical analysis \cite{Butcher1}. A. Connes and D. Kreimer also proved in \cite{ConnesKreimer} that the 
operators $B_+^d$ satisfy the property
\begin{equation}
\label{cocycle}
\Delta\big(B_+^d(t_1\cdots t_n)\big) = 
B_+^d(t_1\cdots t_n)\otimes\un + (\mop{Id}\otimes B^d_+)\circ\Delta(t_1\cdots t_n), 
\end{equation}
for any $t_1,...,t_n\in \Cal T$. This means that $B_+^d$ is a 1-cocycle in the 
Hochschild cohomology of $\HRT$ with values in $\HRT$. 
%%%%%
\section{Simple and contracting arborification}
%%%%%
The Hopf algebra of decorated rooted forests enjoys the following universal property (see e.g. \cite{F}): let $\Cal D$ be a set, let $\Cal H$ be a graded Hopf algebra, and, for any $d\in\Cal D$, let $L^d:\Cal H\to\Cal H$ be a Hochschild one-cocycle, i.e. a linear map such that:
\begin{equation}\label{cocycle}
\Delta\big(L^d(x)\big) = 
L^d(x)\otimes\un_{\Cal H} + (\mop{Id}\otimes L^d)\circ\Delta(x).
\end{equation}
Then there exists a unique Hopf algebra morphism $\Phi:\HRT\to\Cal H$ such that:
\begin{equation}
\Phi_L\circ B_+^d=L^d\circ\Phi_L
\end{equation}
for any $d\in\Cal D$. Now let $V$ be a commutative algebra, let $\big(T(V),\qshu,\Delta\big)$ be the corresponding quasi-shuffle Hopf algebra, let $(e_d)_{d\in \Cal D}$ be a linear basis of $V$, and let $L^d:T(V)\to T(V)$ the right concatenation by $e_d$, defined by:
\begin{equation}
L^d(v_1\ldots v_k):=v_1\ldots v_ke_d.
\end{equation}
One can easily check, due to the particular form of the deconcatenation coproduct, that $L^d$ verifies the one-cocycle condition \eqref{cocycle}. The \textsl {contracting arborification} of the quasi-shuffle Hopf algebra above is the unique Hopf algebra morphism
\begin{equation}
\frak a_V:\HRT\surj \big(T(V),\qshu,\Delta\big)
\end{equation}
such that $\frak a_V\circ B_+^d=L^d\circ\frak a_V$ for any $d\in\Cal D$. It is obviously surjective, since the word $w=e_{d_1}\cdots e_{d_r}$ can be obtained as the image of the ladder $\ell_Y(w)$ with $r$ vertices decorated by $d_1,\ldots d_r$ from top to bottom. This map is invariant under linear base changes. For the shuffle algebra (i.e. when the internal product on $V$ is set to zero), the corresponding Hopf algebra morphism $\frak a_V$ is called \textsl{simple arborification}, and the corresponding section will be denoted by $\ell_X$.\\

Let us apply this construction to multiple zeta values (the base field $k$ being the field ${\mathbb Q}$ of rational numbers): we denote by $\frak a_X$ (resp. $\frak a_Y$) the simple (resp. contracting) arborification from $\HRTX$ onto ${\mathbb Q}\langle X\rangle$ (resp. from $\HRTY$ onto ${\mathbb Q}\langle Y\rangle$). The maps $\zeta_{\sshu}$ and $\zeta_{\sqshu}$ defined in the introduction are characters of the (Hopf) algebras ${\mathbb Q}\langle X\rangle$ and ${\mathbb Q}\langle Y\rangle$ respectively, with values in the algebra ${\mathbb R}[\theta]$. The simple and contracted arborified multiple zeta values are then respectively given by:
\begin{eqnarray}
\zeta_{\sshu}^T:\HRTX&\longrightarrow&{\mathbb R}[\theta]\notag\\
\tau&\longmapsto& \zeta^T_{\sshu}(\tau)=\zeta_{\sshu}\circ\frak a_X(\tau).
\end{eqnarray}
and:
\begin{eqnarray}
\zeta_{\sqshu}^{\ct}:\HRTY&\longrightarrow&{\mathbb R}[\theta]\notag\\
t&\longmapsto& \zeta_{\sqshu}^{\ct}(t)=\zeta_{\sqshu}\circ\frak a_Y(t).
\end{eqnarray}
They extend to any word the maps $\zeta^T$ and $\zeta^{\ct}$ defined in the introduction. Looking back at the examples given there we have:
\begin{equation}
\frak a_Y(\arbrebbdec{n_3}{n_1}{n_2})=y_{n_1}y_{n_2}y_{n_3}+y_{n_2}y_{n_1}y_{n_3}+y_{n_1+n_2}y_{n_3}
\end{equation}
and
\begin{eqnarray*}
\frak a_X(\arbreccbis)&=&2x_0x_0x_1x_1+x_0x_1x_0x_1,\\
\frak a_X(\arbreccter)&=&3x_0x_0x_0x_1.
\end{eqnarray*}
%%%%
\section{Arborification of the map $\frak s$}
%%%%
We are looking for a map $\frak s^T$ which makes the following diagram commutative:
%\begin{equation*}
\diagramme{
\xymatrix{
\Cal H_{\smop{BCK}}^Y\ar@{>>}[d]^{\frak a_Y}\ar[r]^{\frak s^T} & \Cal H_{\smop{BCK}}^X\ar@{>>}[d]^{\frak a_X}\\
\mathbb Q\langle Y\rangle\ar[r]^{\frak s} &\mathbb Q\langle X\rangle
}
}
\noindent An obvious answer to this problem is given by:
$$\frak s^T=\ell_X\circ\frak s\circ \frak a_Y,$$
where $\ell_X$ is the section of $\frak a_X$ described in the previous section. It has the drawback of completely destroying the geometry of trees: indeed, any $Y$-decorated forest is mapped on a linear combination of $X$-decorated ladders. We are then looking for a more natural map with respect to the tree structures, which makes the diagram above commute, or at least the outer square of the diagram below:
\diagramme{
\xymatrix{
\Cal H_{\smop{BCK}}^Y\ar@{>>}[d]^{\frak a_Y}\ar[r]^{\frak s^T}\ar@/_1.5pc/[dd]_{\zeta_{\sqshu}^{T\hskip -6pt -}} & \Cal H_{\smop{BCK}}^X\ar@{>>}[d]^{\frak a_X}\ar@/^2.2pc/[dd]^{\zeta_{\sshu}^{T}} \\
\mathbb Q\langle Y\rangle\ar[r]^{\frak s}\ar[d]^{\zeta_{\sqshu}} &\mathbb Q\langle X\rangle\ar[d]^{\zeta_{\sshu}}\\
\mathbb R[\theta]\ar[r]^{\rho}&\mathbb R[\theta]
}
}
\noindent This interesting problem remains open.
%\end{equation*}


\begin{thebibliography}{abcdsfgh}

	
\bibitem{Br}
	Ch.~Brouder, 
	{\textsl{Runge-Kutta methods and renormalization\/}},
       	Eur.~Phys.~J.~C Part.~Fields 12 (2000) 512--534.
\bibitem{B1}
	F. Brown,
	\textsl{On the decomposition of motivic multiple zeta values},
	Adv. Stud. Pure Math. \textbf{99} (to appear), \texttt{arXiv:1102.1310}.
	
\bibitem{B2}
	F. Brown,
	\textsl{Mixed Tate motives over Spec(${\mathbb Z}$)},
	Duke Math. J. (to appear). \texttt{arXiv:1102.1312}.
	
\bibitem{B3}
	F. Brown,
	\textsl{Depth-graded motivic multiple zeta values},
	\texttt{arXiv:1301.3053} (2013).

\bibitem{Butcher1}
	J.~C.~Butcher,
	{\textsl{An algebraic theory of integration methods}}, 
	 Math.~Comp. 26 (1972) 79--106.
	 
\bibitem{C}
	P. Cartier, \textsl{Fonctions polylogarithmes, nombres polyz\^etas et groupes pro-unipotents}, S\'eminaire Bourbaki No 885, Ast\'erisque \textbf{282}, 137-173 (2002).
	
\bibitem{Cay}
	A.~Cayley,
	{\textsl{On the theory of the analytical forms called trees}}, 
	Phil. Mag. \textbf{13}, 172-176 (1857).

\bibitem{ChaLiv}
    	F.~Chapoton, M.~Livernet,
    	{\textsl{Pre-Lie algebras and the rooted trees operad}},
    	 Internat.~Math.~Res.~Notices 2001 (2001) 395--408.
	 
\bibitem{ConnesKreimer} 
	A.~Connes and D.~Kreimer, 
	\textsl{Hopf Algebras, Renormalization and Noncommutative Geometry.} 
	Comm.\ Math.\ Phys.\ {\bf 199} (1998) 203-242.
	
\bibitem{DG}P. Deligne, A. Goncharov,
	\textsl{Groupes fondamentaux motiviques de Tate mixtes},
	Ann. Sci. Ec. Norm. Sup. (4) \textbf{38} No1, 1-56 (2005).
	 
\bibitem{Dur}
	A. D\" ur,
	\textsl{M\"obius functions, incidence algebras and power series representations},
	Lect. Notes math. \textbf{1202}, Springer (1986).

\bibitem{E1}
	J. Ecalle,
	\textsl{Les fonctions r\'esurgentes} Vol. 1, Publications Math\'ematiques d'Orsay (1981). Available at http://portail.mathdoc.fr/PMO/feuilleter.php?id=PMO\_ 1981.
	
\bibitem{E2}
	J. Ecalle,
	\textsl{Les fonctions r\'esurgentes} Vol. 2, Publications Math\'ematiques d'Orsay (1981). Available at http://portail.mathdoc.fr/PMO/feuilleter.php?id=PMO\_ 1981.
	
\bibitem{E3}
	J. Ecalle,
	\textsl{Les fonctions r\'esurgentes} Vol. 3, Publications Math\'ematiques d'Orsay (1985). Available at http://portail.mathdoc.fr/PMO/feuilleter.php?id=PMO\_ 1985.	
	
\bibitem{E92}
	J. Ecalle,
	\textsl{Singularit\'es non abordables par la g\'eom\'etrie},
	Ann. Inst. Fourier \textbf{42}, No 1-2, 73-164 (1992).
	
\bibitem{EV}
	J. Ecalle, B. Vallet,
	\textsl{The arborification-coarborification transform: analytic, combinatorial, and algebraic aspects},
	Ann. Fac. Sci. Toulouse \textbf{XIII}, No 4, 575-657 (2004).

\bibitem{FM}
	F. Fauvet, F. Menous,
	\textsl{Ecalle's arborification-coarborification transforms and the Connes-Kreimer Hopf algebra},
	Ann. Sci. Ec. Norm. Sup., to appear.
	\texttt{arXiv:1212.4740} (2012).
	
\bibitem{F}
	L. Foissy,
	\textsl{Les alg\`ebres de Hopf des arbres enracin\'es d\'ecor\'es I,II},
	Bull. Sci. Math. \textbf{126} (2002), 193-239 and 249-288.

	
\bibitem{H}
	M. E. Hoffman,
	\textsl{Multiple harmonic series},
	Pacific J. Math. \textbf{152}, 275-290 (1992).
	
\bibitem{H2}
	M. E. Hoffman,
	\textsl{Quasi-shuffle products},
	J. Algebraic Combin. \textbf{11}, 49-68 (2000).
	
\bibitem{IKZ}
	K. Ihara, M. Kaneko, D. Zagier,
	\textsl{Derivation and double shuffle relations for multiple zeta values},
	Comp. Math. \textbf{142} No2, 307-338 (2004).
	
\bibitem{Kreimer} 
	D.~Kreimer, 
	\textsl{On the Hopf algebra structure of perturbative quantum field theories}, 
	Adv. Theor. Math. Phys.{\bf 2}, 303--334 (1998).
	
\bibitem{Mu}
	A. Murua, 
	\textsl{The Hopf algebra of rooted trees, free Lie algebras, and Lie series}, 
	Found. Computational Math. \textbf{6}, 387-426 (2006).

\bibitem{R}G. Racinet,
	\textsl{Doubles m\'elanges des polylogarithmes multiples aux racines de l'unit\'e},
	Publ. Math. IHES \textbf{95}, 185-231 (2002).
	
\bibitem{T}T. Terasoma,
	\textsl{Mixed Tate motives and multiple zeta values},
	Invent. Math. \textbf{149} No2, 339-369 (2002).
	
\bibitem{W}
	M. Waldschmidt, \textsl{Valeurs z\^eta multiples, une introduction},
	J. Th\'eor. Nombres de Bordeaux \textbf{12}, 581-595 (2000).
	
\bibitem{Y14}S. Yamamoto,
	Multiple zeta-star values and multiple integrals, preprint, \texttt{arXiv:1405.6499} (2014).
	 
\bibitem{Z}
	D. Zagier, \textsl{Values of zeta functions and their applications}, Proc. First European Congress 	of Mathematics, Vol. 2, 497-512, Birkh\" auser, Boston (1994).

\end{thebibliography}
\end{document}